\documentclass [12pt]{article}
\newfont{\bbb} {msbm10}
\newcommand{\Bbb}[1]{\mbox{\bbb#1}}
\newcommand{\R}{\Bbb{R}}

\newcommand{\ra}{\rightarrow}
\newcommand{\M}{{\cal{M}}}
\newcommand{\bM}{\bar{{\cal{M}}}}
\newcommand{\hM}{\hat{{\cal{M}}}}

\newcommand{\cM}{{\cal{M}}}
\newcommand{\cC}{{\cal{C}}}

\usepackage[all]{xy}

\begin{document}

\title{Exotic Structures and the Limitations of Certain Analytic Methods in Geometry}
\author{F. T. Farrell and P. Ontaneda\thanks{The first author was
partially supported by a NSF grant.
The second author was supported in part 
by a research  grant from CAPES, Brazil.}}
\date{}

\maketitle

In this survey we review some results concerning negatively curved exotic 
structures ($DIFF$ and $PL$) and its (unexpected) implications on the limitations of some analytic
methods in geometry. Among these methods are the harmonic map method and the
Ricci flow method. \\

First in section 1 we mention certain results about the rigidity of negatively curved manifolds.
In section 2  and 3 we survey some results concerning the limitations of the harmonic map technique
and the natural map technique for negatively curved manifolds. Finally, in section 4, we 
mention some limitations of the Ricci flow method for pinched negatively curved manifolds.\\

We are grateful to J-F. Lafont and R. Spatzier for the useful information they provided to us
and to E. Gasparim for suggesting some improvements in the text. We are also grateful to the 
referee for pointing out certain inaccuracies. 
\vspace{.6in}

{\large {\bf 1. Negative curvature and rigidity.}}\\

We begin with a basic question in geometry and topology:\\

{\it When are two homotopy equivalent manifolds diffeomorphic, $PL$ homeomorphic or homeomorphic?}\\

If both manifolds are closed, hyperbolic and of dimension greater than 2, Mostow's Rigidity Theorem \cite{Mos}
says that they are isometric, in particular diffeomorphic. When both manifolds have strictly negative curvature, 
results of Eells and Sampson \cite{ES}, Hartman \cite{Har} and Al'ber \cite{A'} show that if 
$f:M_{1}\rightarrow M_{2}$ is a homotopy equivalence then it is homotopic to a unique harmonic map
(see also next section). 
Lawson and Yau conjectured that this harmonic map is always
a diffeomorphism (see problem 12 of a list of problems presented by Yau in \cite{Yau}). 
Farrell and Jones \cite{FJ1} gave counterexamples to
this conjecture by proving the following. \\

\noindent {\bf Theorem 1. \cite{FJ1}} {\it If $M$ is a real hyperbolic manifold
and $\Sigma$ is an exotic sphere, then given $\epsilon>0$, $M$ has a finite
covering $\tilde{M}$ such that  the connected  sum $\tilde{M}\#\Sigma$ is not diffeomorphic to
$\tilde{M}$ and  admits a Riemannian metric with all sectional curvatures in the interval $(-1-\epsilon,-1+\epsilon)$.}\\

Since there are no exotic spheres in dimensions $<7$  this does not give
counterexamples to Lawson-Yau conjecture in dimensions less than 7. 
(Also note that, for example, there are no exotic 12-dimensional spheres.) 
Moreover, since the $DIFF$ category is equivalent to the
$PL$ category in dimensions less than 7, changing the differentiable structure is equivalent to changing the $PL$
structure. The Theorem above was generalized by Ontaneda in \cite{O} to dimension 6, by changing the $PL$ structure, and this
result was extended in \cite{FJO1} to all dimensions greater that five:\\

\noindent {\bf  Theorem 2. \cite{FJO1}}  {\it For every $n>5$, there are closed real hyperbolic $n$-manifolds $M$  
such that the following holds. Given $\epsilon>0$, $M$ has a finite cover $\tilde{M}$ that supports
an exotic smoothable $PL$ structure that  admits a Riemannian metric with sectional curvatures in the
interval $(-1-\epsilon,-1+\epsilon)$.}\\

This result gives counterexamples to the Lawson-Yau conjecture in all dimensions $>5$ since
Whitehead showed that a smooth manifold has a unique $PL$ structure. (Recall that two smooth
manifolds are $PL$ equivalent if and only if there is a simplicial complex which smoothly triangulates
both manifolds.)\\

The hyperbolic manifolds mentioned in the Theorem above are obtained using methods of
Millson and Raghunathan \cite{MR}, which were based on a earlier work of Millson \cite{Mill}. These methods provide
a large class of examples of hyperbolic manifolds with many non-vanishing cohomology classes.
\vspace{.6in}

{\large {\bf 2. Negative curvature and harmonic maps.}}\\

Let $M$ and $N$ be two compact Riemannian manifolds. Recall that the
energy of a map $f: M\ra N$ is defined to be $\frac{1}{2}\int_{M}|df|^2$.
A harmonic map is a map which is a  critical point of this functional.
It satisfies the equation $\tau (f)=0$, where $\tau(f)$ is the tension field of $f$,
(see for example \cite{EL2}, p.14.)\\

Part of the interest in harmonic maps comes from the fact that they are very successful in 
proving rigidity (and superrigidity) results for non-positively curved Riemannian manifolds.  
We can mention for
example results of Siu \cite{Siu}, Sampson \cite{Sa2}, Hern{\'a}ndez \cite{Her},
Corlette \cite{Cor}, Gromov and
Schoen \cite{GS}, Jost and Yau \cite{JY}, and Mok, Sui and Yeung \cite{MSY}.  All of which are
based on the pioneering existence Theorem of Eells and Sampson \cite{ES} and the uniqueness
Theorem of Hartman \cite{Har} and Al'ber \cite{A'}. Eells and Sampson proved that 
given any smooth map $k_0 : M\ra N$ between Riemannian manifolds, 
the heat flow equation, that is, the PDE initial value problem
\begin{equation}
\frac{\partial k_t}{\partial t} = \tau(k_t),\ \ \ k_t\big|_{t=0} = k_0
\end{equation}
has a unique solution $k_t$ (for all $t\geq 0$) and that
$\lim_{t\to\infty}k_t = k$; cf.\cite{EL2}, pp. 22-24.  
Here $N$ is assumed to have non-positive curvature.
Note that $k_t$ is a homotopy between $k_0$ and $k$ and that $k$ is a harmonic map.
Also, if in addition, $M$ has negative curvature the results of Eells and Sampson together with the
results of Hartman \cite{Har} and Al'ber \cite{A'} show that there is a unique harmonic map
homotopic to $k_0$.
\vspace{.4in}

\noindent {\bf 2.1 Lawson-Yau conjecture.}\\

Let $f:M\rightarrow N$ be a homotopy equivalence between negatively curved manifolds.
As already mentioned in section 1, Lawson and Yau conjectured that the unique harmonic map
$\phi : M\rightarrow N$ 
homotopic to $f$ is a diffeomorphism. Theorems 1 and 2 proved that this conjecture is
false in dimensions $>5$. That is, for every dimension $>5$ there are harmonic homotopy
equivalences $f:M\rightarrow N$ which are not diffeomorphisms.
Theorems 1 and 2 already place some limitations to the harmonic maps technique.
But there  remained the question whether a ``topological" Lawson-Yau conjecture could hold:
\vspace{.2in}

\noindent {\bf (*)} {\it Let $\phi : M\rightarrow N$ be a harmonic homotopy equivalence between closed
negatively curved manifolds. Is $\phi$ a homeomorphism?}.
\vspace{.2in}

A positive answer to this conjecture would give an analytic proof of ``Borel's Conjecture" for 
closed negatively curved manifolds:
\vspace{.2in}

\noindent {\bf Borel's Conjecture.} {\it Let $M$ and $N$ be homotopy equivalent  closed aspherical
manifolds. Then $M$ and $N$ are homeomorphic.}
\vspace{.2in}

Borel's conjecture has been verified in \cite{FJ2} when one of the manifolds is non-positively 
curved and dimensions $\neq$ 3,4. The proof uses sophisticated topological methods.
On the other hand, a negative answer to {\bf (*)}, would imply that this last result (the proof 
of Borel's conjecture for closed non-positively curved manifolds in \cite{FJ2})
cannot be obtained, at least directly, using the harmonic maps technique.\\

\noindent {\bf Remark.} Conjecture {\bf (*)} was studied in \cite{FJ3} and some partial (negative) results were given.\\

But the topological Lawson-Yau conjecture {\bf (*)} is also false:\\

\noindent {\bf Theorem 3. \cite{FJO2}} {\it In every dimension $n\geq 6$, there is a pair of closed negatively curved
 manifolds $M^{n}$ and $N^{n}$ and a harmonic homotopy equivalence $\phi : M^{n}\rightarrow
N^{n}$, which is not one-to-one.}\\

Actually, we can prove a little more:\\

\noindent {\bf Theorem 4. \cite{FJO2}} {\it In every dimension $n\geq 6$, there is a pair of closed negatively curved
 manifolds $M^{n}$ and $N^{n}$ such that the following holds. For any  homotopy equivalence $f:M^n\ra N^n$, the
 unique harmonic map $\phi : M^{n}\rightarrow
N^{n}$ homotopic to $f$ is not one-to-one.}\\

This Theorem can be directly deduced from Theorem 2 and the $C^\infty -$ Hauptvermutung of Scharlemann
and Siebenmann \cite{SC1}. We reproduce this short deduction here since it shows, unexpectedly,  how
the theory of $PL$ manifolds interweaves with the theory of harmonic maps.\\
 
\noindent {\bf Proof.} By Theorem 2 we have\\

\noindent {\bf (2.1.1)} {\it In every dimension $n\geq 6$, there is a pair of non-PL-equivalent closed negatively curved
 manifolds $M^{n}$ and $N^{n}$ with $\pi_{1}( M^{n})$ isomorphic to $\pi_{1}( N^{n})$.} \\

Let $M^{n}$ and $N^{n}$ be a pair of manifolds satisfying {\bf (2.1.1)}, and let $\phi : M^{n}\rightarrow N^{n}$ be 
the unique harmonic map realizing the isomorphism 
$\pi_{1}(M^{n})\rightarrow \pi_{1}(N^{n})$ induced by the homotopy equivalence $f$. The Theorem now follows by 
just applying the following result of M. Scharlemann and L. Siebenmann \cite{SC1}\\

\noindent {\bf (2.1.2)} {\it Smoothly homeomorphic closed manifolds of dimension $\geq$ 6 are PL-homeomorphic.}\\

\noindent {\bf Remark.} Smooth homeomorphisms are not necessarily diffeomorphisms. A simple example is given by
the smooth homeomorphism $f:\R\rightarrow\R$, $f(x)=x^3$.\\

Thus, the harmonic map $\phi$ cannot be a homeomorphism because $M^{n}$ and $N^{n}$ satisfy
{\bf (2.1.1)}. This proves the Theorem.\\

Recall that Poincar{\'e} Conjecture in low dimensional topology asserts that the only simply 
connected closed 3-dimensional
manifold is the 3-sphere (up to homeomorphism), or, equivalently, that every homotopy 3-sphere is homeomorphic to $S^{3}$. 
Now, by using another result of M. Scharlemann (and assuming Poincare's conjecture) we can get a little more:\\

\noindent {\bf Theorem 5. \cite{FJO2}} {\it Assume that every homotopy 3-sphere is homeomorphic to $S^{3}$. 
Then in every dimension $n\geq 6$, there is a pair of closed negatively curved
manifolds $M^{n}$ and $N^{n}$ and a harmonic homotopy equivalence $\phi : M^{n}\rightarrow
N^{n}$, which is not cellular}\\

See \cite{E} for a discussion of cellular maps (which are called cell like
maps in that article). Siebenmann \cite{Si} showed that a continuous map $f : X\to Y$
between a pair of closed manifolds of dimension $\geq 5$ is cellular if and only if it is the
limit of homeomorphisms.\\

As with Theorem 3, Theorem 5 is a direct consequence of little more general one:\\

\noindent {\bf Theorem 6. \cite{FJO2}} {\it Assume that every  homotopy 3-sphere is homeomorphic to $S^{3}$. 
Then in every dimension $n\geq 6$, there is a pair of closed negatively curved
manifolds $M^{n}$ and $N^{n}$ such that the following holds. For any  homotopy equivalence $f:M^n\ra N^n$, the
unique harmonic map $\phi : M^{n}\rightarrow N^{n}$ homotopic to $f$ is not cellular, i.e. it is not 
the uniform limit of homeomorphisms.

Consequently the maps $k_t$ and $l_t$ in
the heat flow of  $f= k_0$ to $k = k_\infty$ and of $g= l_0$ to $l =
l_\infty$ are not one-to-one for all $t$ sufficiently large.
Here $g$ is a homotopy inverse to $f$.}\\

The proof is the same as the one in Theorem 4, but now we use {(2.1.2)} together with (see \cite{SC2}):\\

\noindent {\it Assume that every homotopy 3-sphere is homeomorphic to $S^{3}$. Then any smooth
cellular map $\phi :M^{n}\rightarrow N^{n}$ of smooth closed n-manifolds
(where $n\geq 6$) is smoothly homotopic,
through cellular maps, to a smooth homeomorphism}.\\

\noindent {\bf Remark.} In all the Theorems of this subsection we can assume that one of the manifolds
is hyperbolic. This follows from Theorem 2.
\vspace{.4in}

\noindent {\bf 2.2. Yau's problem 111.}\\

Let $f:M_1\rightarrow M_2$ be a homotopy equivalence between negatively curved manifolds
and let $h:M_1\rightarrow M_2$ be the unique harmonic map homotopic to $f$.
In the examples provided by the Theorems above, the main obstruction to $h$ being a diffeomorphism 
or a homeomorphism is that $M_1$ and $M_2$ are not $PL$ equivalent, even though they
are homotopy equivalent (in fact homeomorphic). We may ask then what happens if this
obstruction vanishes, that is, if $M_1$ and $M_2$ are diffeomorphic. Can the harmonic
map technique be applied in this context to obtain diffeomorphisms or, at least, 
homeomorphisms? Or, equivalently, if we flow a diffeomorphism (using the heat flow),
will the limit be also a diffeomorphism or a homeomorphism?
This is considered in  Problem 111 of the list compiled by S.-T.  
Yau in  \cite{Yau}.  Here is a restatement of this problem.\\

\noindent {\bf Problem 111 of \cite{Yau}}.  Let $f:M_1\rightarrow M_2$ be a
diffeomorphism between two compact manifolds with negative curvature.  If
$h:M_1\rightarrow M_2$ is the unique harmonic map which is homotopic to $f$, is $h$ a
homeomorphism?, or equivalently, is $h$ one-to-one?\\

(This problem had been reposed in \cite{Sta} as Grand Challenge Problem 3.6.)  The
answer to the problem was proved to be yes when dim$M_1=2$ by Schoen-Yau
\cite{SY} and Sampson  \cite{Sa}.  But it was proved by Farrell, Ontaneda and Raghunathan
\cite{FOR} that the answer to this question is negative.\\

\noindent {\bf Theorem 7. \cite{FOR}} {\it For every integer $n\geq 6$, there is a diffeomorphism $f:M_1
\rightarrow M_2$ between a pair of closed negatively curved $n$-dimensional
Riemannian manifolds such that the unique harmonic map $h:M_1 \rightarrow
M_2$ homotopic to $f$ is not one-to-one.}  \\

\noindent {\bf Addendum.}  {\it In the Main Theorem, either $M_1$ or $M_2$ can be chosen to
be a real hyperbolic manifold and the other chosen to have its sectional
curvatures pinched within $\varepsilon$ of $-1$; where $\varepsilon$ is any
preassigned positive number.}\\

Hence the negative answer given by this Theorem to Problem 111 places more limits to the
applicability of the harmonic map technique to rigidity questions. \\

Theorem 7 evolves from Theorems 1-6 above and follows from Theorem 8 below.\\

\noindent {\bf Theorem 8. \cite{FOR}} {\it Given an integer $n\geq 6$ and a positive real number
$\varepsilon$, there exists a $n$-dimensional closed connected orientable
(real) hyperbolic manifold $M$ and a homeomorphism $g:\cM \rightarrow M$ with
the following properties:}

\begin{enumerate}
\item[{1.}]  $\cM$ {\it is a negatively curved Riemannian manifold whose sectional
curvatures are all in the interval} $(-1-\varepsilon, -1+\varepsilon)$.

\item[{2.}] {\it $M$ and $\cM$ are not $PL$ homeomorphic.  }

\item[{3.}] {\it There is a connected 2-sheeted covering space ${\tilde {M}}
\rightarrow M$ such that ${\tilde {g}} : {\tilde {\cM}} \rightarrow {\tilde
{M}}$ is homotopic to a diffeomorphism.}
\end{enumerate}

\noindent {\bf Remark.}  In property 3, ${\tilde {\cM}} \rightarrow \cM$ denotes
the pullback of the covering space ${\tilde {M}} \rightarrow M$ via $g$, and
${\tilde {g}}$ is the induced homeomorphism making the diagram

$$
\begin{array}{ccc} {\tilde {\cM}} & \stackrel{\tilde {g}}  \rightarrow & {\tilde
{M}} \\
\downarrow & &\downarrow \\
\cM & \stackrel{g} \rightarrow & M  
 \end{array}
$$
into a Cartesian square.  Also, ${\tilde {M}}$ and ${\tilde {\cM}}$ are
given the differential structure and Riemannian metric induced by ${\tilde
{M}} \rightarrow M$ and ${\tilde {\cM}} \rightarrow \cM$, respectively.\\

The key ingredient in the proof of Theorem 8 is the existence of closed real 
hyperbolic manifolds with interesting cup product properties.  
Such manifolds are constructed in section 2 of \cite{FOR}.\\

\noindent {\bf Proof of Theorem 7 assuming Theorem 8.} Let $g:\cM \rightarrow M$ be the
homeomorphism given by Theorem 8 relative to $n$ and $\varepsilon$.  Set
$M_1 = {\tilde {\cM}}, M_2 = {\tilde {M}}$ and let $f:M_1 \rightarrow M_2$ be a
diffeomorphism homotopic to ${\tilde {g}}: {\tilde {\cM}} \rightarrow {\tilde
{M}}$ which exists by property 3 of Theorem 8.  Let $k:\cM \rightarrow M$ be
the unique harmonic map homotopic to $g$  given by the fundamental existence
result of Eells and Sampson \cite{ES} and uniqueness by Hartmann \cite{Har} and Al'ber 
\cite{A'}.  Lifting this homotopy to the covering spaces ${\tilde {\cM}}, {\tilde
{M}}$ gives a smooth map 
$$
{\tilde {k}}: {\tilde {\cM}} \rightarrow {\tilde {M}}$$
covering $k$ and homotopic to ${\tilde {g}}$.  Note that ${\tilde {k}}$ is
also a harmonic map as is easily deduced from \cite{EL1}, 2.20 and 2.32.
Consequently, ${\tilde {k}}$ is the harmonic map $h:M_1 \rightarrow M_2$
mentioned in the statement of Theorem 7.  Also note that if ${\tilde {k}}$
is univalent, then so is $k$.  Hence it suffices to show that $k$ is   {\it
not} univalent.  Since $k$ is smooth, $k$ univalent would mean that 
$$
k:\cM \rightarrow
M$$
is a $C^\infty$-homeomorphism and hence $M$ and $\cM$ are $PL$-homeomorphic
by the $C^\infty$-Hauptvermutung proved by  Scharlemann and Siebenmann \cite{SC1}.
And this would contradict property 2 of Theorem 8;  consequently, $k$ and
hence also $h$ are {\it not} univalent.  This   proves the  Theorem 7 and the
part of the Addendum where $M_2$ is real hyperbolic.

To do the case where $M_1$ is real hyperbolic; set $M_1 = {\tilde {M}}, M_2
= {\tilde {\cM}}$ and let $f$ be a diffeomorphism homotopic to ${\tilde
{g}}^{-1}$.  The rest of the argument is  a before. This concludes the deduction
of Theorem 7 from Theorem 8.\\

Note that, as in Theorems 3-6, crucial use is made here of the Scharlemann-Siebenmann
$C^\infty$-Hauptvermutung \cite{SC1}. \\

Hence the idea of the proof of Theorem 7 can be paraphrased in the following few words.
Take a homotopy equivalence $f:M_1\ra M_2$ between homeomorphic negatively curved
manifolds, with $M_1$ not $PL$-homeomorphic to $M_2$. Theorem 2 grants the existence
of such objects in every dimension $>5$. Then, as shown in Theorem 4, the unique harmonic
map $h:M_1\ra M_2$ homotopic to $f$ cannot be one-to-one. Suppose that after taking
some finite cover ${\tilde{f}}:{\tilde{M}}_1\ra{\tilde{M}}_2$ $\tilde f$ becomes
homotopic to a diffeomorphism. Let $k$ be the unique harmonic map homotopic
to $\tilde f$. Then $k$ is homotopic to a diffeomorphism but $k$ is not one-to-one
since (even though the $PL$ obstruction now vanishes) the damage is already done:
$k=\tilde h$. The existence of manifolds admitting
such finite covers (in fact double covers) is granted by Theorem 8.
\vspace{.3in}

\noindent {\bf 2.3. Cellular harmonic maps.}\\

Since a harmonic map (between closed negatively curved manifolds)
homotopic to a diffeomorphism is not necessarily a homeomorphism we can
ask a deeper question: suppose now that the harmonic map can be approximated
by homeomorphisms (or even diffeomorphisms), that is, the harmonic map
is cellular. Does this imply that the harmonic map is a diffeomorphism?
The following Theorem shows that the answer to this question is also
negative, showing even more limitations to the harmonic map technique:\\

\noindent {\bf Theorem 9. \cite{FO1}} {\it For every integer $m>10$, there is a 
harmonic cellular map  $h:M_{1}\ra M_{2}$, between a pair
of closed negatively curved $m$-dimensional Riemannian
manifolds, which is not a diffeomorphism.}\\

\noindent {\bf Addendum.} {\it The map $h$ in Theorem 9 can be approximated 
by diffeomorphisms. Also, either $M_{1}$ or
$M_{2}$ can be chosen to be a real hyperbolic manifold and
the other chosen to have its sectional curvatures pinched within 
$\epsilon$ of -1; where $\epsilon$ is any preassigned positive
number.}\\

We conjecture that this can be improved to all dimensions $\geq 6$.
We do not know whether the harmonic map $h$ in the statement of Theorem 9 can ever be
a homeomorphism. \\

Theorem 9 follows from the next Theorem,  which is of independent interest:\\

\noindent {\bf Theorem 10. \cite{FO1}} {\it For every integer $m>10$, and
$\epsilon >0$, there are an $m$-dimensional closed orientable smooth
manifold $\M$, and a $C^{\infty}$ family of Riemannian metrics $\mu_{s}$, 
on $\M$, $s\in [0,1]$, such that:}

\begin{enumerate}
\item[{(i)}] {\it $\mu_{1}$ is hyperbolic. }

\item[{(ii)}] {\it The sectional curvatures of $\mu_{s}$, $s\in [0,1]$, 
are all in interval $(-1-\epsilon, -1+\epsilon)$. }

\item[{(iii)}] {\it The maps $k$ and $l$ are both not univalent (i.e.\ not one-to-one) where $k :
(\M,\mu_0) \to (\M,\mu_1)$ and $l : (\M,\mu_1) \to (\M,\mu_0)$ are the unique harmonic maps
homotopic to \mbox{\rm id}$_\M$.} 
\end{enumerate}

The derivation of Theorem 9 from Theorem 10 uses the
continuous dependence (in the $C^\infty$-topology) of the harmonic map homotopic to a homotopy
equivalence $f : (M,\mu_M)\to (N,\mu_N)$ on the negatively curved Riemannian metrics $\mu_M$
and $\mu_N$.  This dependence was proved by Sampson \cite{Sa}, Schoen and Yau \cite{SY}, and
Eells and Lemaire \cite{EL2}. To derive Theorem 9 from Theorem 10 let $k_{t}:
(\M,\mu_0) \to (\M,\mu_t)$ be the unique harmonic map homotopic to id. Then $k_0 =$id
and $k_1 =k$, which is not one-to-one. Since the space of diffeomorphisms is open in the $C^k$ 
topology ($k\geq 1$) it follows that there is a minimal $t_0 >0$ such that $k_{t_0}$ is not a
diffeomorphism and $k_{t_0}$ can be approximated by the diffeomorphisms $k_{t}$, $t <t_0$.\\

Likewise, as with Theorems 3 and 4, Scharlemann's result \cite{SC2} also implies
another curious relationship between Poincar{\'e} Conjecture in low
dimensional topology and the existence of a certain type of harmonic map $k :M\to N$ between
high dimensional (i.e.\ dim $M >10$) closed negatively curved Riemannian manifolds.   
If Poincar{\'e} Conjecture holds, then there exists a harmonic map $k$ which is
homotopic to a diffeomorphism but cannot be approximated by homeomorphisms; i.e. is not a
cellular map.  Explicitly, we have the following addendum to Theorem 10:\\

\noindent {\bf Addendum to Theorem 10.}  {\it Assuming that the Poincar{\'e} Conjecture is true, then 
the harmonic maps $k$ and
$l$ (of Theorem 10) are not cellular.  And consequently the maps $k_t$ and $l_t$ in
the heat flow of {\rm id} $= k_0$ to $k = k_\infty$ and of {\rm id} $= l_0$ to $l =
l_\infty$ are not univalent for all $t$ sufficiently large.}\\

The key to the proof of Theorem 10 is the following important result,
which is also used in the proofs of the results of section 4 that show some limitations of the Ricci
flow method:\\

\noindent {\bf Theorem 11 \cite{FO1}.} {\it
Given an integer $m>10$ and a positive number $\epsilon$, there
exist a $m$-dimensional closed orientable real hyperbolic manifold
$M$ and a smooth manifold $\M$ with the following properties:}

\begin{enumerate}
\item[(i)] {\it $M$ is  homeomorphic to $\M$.}

\item[(ii)] {\it $M$ is not PL homeomorphic to $\M$.}

\item[(iii)]{\it $\M$ admits a Riemannian metric $\mu$, whose sectional 
curvatures are all in the interval $(-1-\epsilon , -1+\epsilon )$.}

\item[(iv)] {\it There is a finite sheeted cover  $p:\bM\ra\M$ and a 
one-parameter $C^{\infty}$ family of Riemannian metrics $\mu_{s}$, 
on $\bM$, $s\in [0,1]$, such that $\mu_{0}=p^{*}\mu$ and $\mu_{1}$ is
hyperbolic. The sectional curvatures of $\mu_{s}$, $s\in [0,1]$, 
are all in the interval $(-1-\epsilon, -1+\epsilon)$.}
\end{enumerate}

The proof of Theorem 10 assuming Theorem 11 resembles the proof of Theorem
7 (assuming Theorem 8) given before. Again, crucial use is made of the
$C^{\infty}$-Hauptvermutung of Scharlemann-Siebenmann \cite{SC1}.\\

We outline the proof of Theorem 11. By Theorem 8 there is a pair of homeomorphic but not PL 
homeomorphic closed negatively curved Riemannian manifolds $M$ and $\M$ satisfying:
\begin{enumerate}
\item  $M$ is real hyperbolic.
\item  $\M$ has a 2-sheeted cover $q : \hM \to \M$ where
$\hM$ admits a real hyperbolic metric $\nu$.
\end{enumerate}
Let $\mu$ be a given negatively curved Riemannian metric on $\M$ and $q^*(\mu)$ be the induced
Riemannian metric on $\hM$.  We would like to find a 1-parameter family of negatively curved
Riemannian metrics connecting $q^*(\mu)$ to $\nu$.  But we don't know how to do this.  In fact
this is in general an open problem \cite[Question 7.1]{BK}.  However by passing to a
large finite sheeted cover $r : \bM \to \hM$, we are able to connect $(q\circ r)^*(\mu)$ to
the real hyperbolic metric $r^*(\nu)$ by a 1-parameter family of negatively curved Riemannian
metrics; this is essentially the content of Theorem 11 in which $p=q\circ r$.  To
accomplish this, several results about smooth pseudo-isotopies are used; in particular, the
main result of \cite{FJ} concerning the space of stable topological pseudo-isotopies of real
hyperbolic manifolds together with the comparison between the spaces of stable smooth and
stable topological pseudo-isotopies contained in \cite{BL} and \cite{Hat}.  And finally we need
Igusa's fundamental result \cite{I} comparing the spaces of pseudo-isotopies and stable
pseudo-isotopies.  We need that dim $M > 10$ in order to invoke Igusa's result. 
\vspace{.6in}

{\large {\bf 3. Natural maps and negative curvature.}}\\

It was pointed out to us by M. Varisco \cite{V} that the limitations of the harmonic
map technique obtained by the results of section 2 can also be applied
to the {\it natural maps} defined by G. Besson, G. Courtois and S. Gallot
\cite{BCG}.\\

Given a homotopy equivalence  $f:M\ra N$  between closed negatively
curved manifolds G. Besson, G. Courtois and S. Gallot \cite{BCG} defined
the {\it natural map}  $f^*:M\ra N$ associated to $f$. The map $f^*$ has
many interesting geometric and dynamic properties. Like harmonic maps
they are also useful for proving rigidity results. For example
Mostow's Rigidity Theorem for hyperbolic manifolds \cite{Mos} 
follows from the following Theorem:\\

\noindent {\bf Theorem. \cite{BCG}} {\it Let $f:M\ra N$  be a homotopy
equivalence between closed negatively
curved locally symmetric spaces of dimension $\geq 3$. Then 
(possibly after rescaling) the natural map $f^*:M\ra N$  is an isometry.}\\

But we are interested in the following properties of natural maps:\\
\begin{enumerate}
\item[{1.}] $f^*$ is at least $C^1$ (\cite{BCG}, p. 635).

\item[{2.}]  $f^*$ is homotopic to $f$ (\cite{BCG}, p.634).

\item[{3.}] If $f$ is homotopic to $g$ then $f^*=g^*$.

\item[{4.}] If $\bar{f}: \bar{M}\ra\bar{N}$ is a finite cover of $f:M\ra N$ then $(\bar{f})^*=\overline{f^*}$.
%\item[{5.}] $f^*$ depends continuously on $f$ and on the metrics of $M$ and $N$.
%Here we consider $f^*$ and $f$ varying in the $C^1$ topology and the metrics varying in the $C^2$ topology.
\end{enumerate}

Property 3 holds because $\partial {\tilde{f}}=\partial {\tilde{g}}:
\partial {\tilde{M}}\ra \partial {\tilde{N}}$ (see \cite{BCG}, pp.633-634).
Property 4 follows directly from the definition of natural maps. \\

It may be argued that natural maps are, in some sense, better than
harmonic maps. But, as observed by M. Varisco, property 2 above
implies that Theorem 4 (with its addendum) also holds for natural maps:\\

\noindent {\bf Theorem 12.} {\it In every dimension $n\geq 6$, there is a pair of closed negatively curved
 manifolds $M^{n}$ and $N^{n}$ such that the following holds. For any  homotopy equivalence $f:M^n\ra N^n$, the
 the natural map $f^* : M^{n}\rightarrow
N^{n}$ is not one-to-one.}\\

Also, a version of Theorem 6 holds for natural maps. Note also that properties 1,2,3,4 imply 
a version of Theorem 7 for natural maps:\\

\noindent {\bf Theorem 13.} {\it For every integer $n\geq 6$, there is a diffeomorphism $f:M_1
\rightarrow M_2$ between a pair of closed negatively curved $n$-dimensional
Riemannian manifolds such that the natural map $f^*:M_1 \rightarrow
M_2$ is not one-to-one.}  \\

(It is the first time that the statements of Theorems 12 and 13 appear in print.
The proofs are similar to the proofs of Theorems 4 and 6, respectively.)\\

{\bf Remark.} We do not know whether versions of Theorems 9 and 10 hold for natural maps.
In particular, we do not know if there are natural maps that can be approximated
by diffeomorphisms but are not diffeomorphisms. This is an interesting question.
To have  versions of Theorems 9 and 10 hold for natural maps
we would need to show that  $f^*$ depends continuously on the metrics of $M$ and $N$,
where we consider $f^*$ varying in the $C^1$ topology and the metrics varying in the $C^2$ topology.
(In fact, in our examples the metrics vary in the $C^{\infty}$ topology.)
One way one may try to verify this continuous dependence would be to use the fact that
the natural map is defined implicitly by the equation $G(F(y),y)=0$
in p.636 of \cite{BCG} ($F$ is the natural map in this equation).
The entropy of one of the metrics and the Busemann functions appear in the definition of
the function $G$.  Note that the perturbations of the entropy (with respect to the metric) have some regularity 
(see \cite{KPW}). We could not find a reference for the regularity of the perturbations of the
Busemann functions (with respect to the metric) but the proof of the smoothness of the Busemann
functions (for universal covers of closed smooth negatively curved manifolds) that appears
in \cite{Sh} might be useful.\\

All this can be generalized. The following definition tries to formalize any
process (analytic or otherwise) that assigns to every continuous map between
closed negatively curved manifolds a {\it special} map. For manifolds $M,N$, we denote the space of
continuous maps $M\ra N$ by $\cC (M,N)$.\\

\noindent {\bf Definitions.} A {\it special correspondence} $\Psi $ for closed
negatively curved manifolds is just a family of maps $\Psi_{M,N}:\cC(M,N)\ra \cC(M,N)$, for
each pair of closed negatively curved manifolds $M, $ $N$.
Note that $\Psi_{M,N}$ depends on the metrics of $M$ and $N$.  For $f:M\ra N$, we say that
$\Psi f$ is the $\Psi$-special map associated to $f$. We say that $\Psi$ is $C^k$ if
$\Psi f$ is $C^k$, for every $f$. We say that $\Psi$ is a homotopy special
correspondence if $\Psi f$ is homotopic to $f$, for every $f$, and $\Psi f=\Psi\, g$
for every $f$ homotopic to $g$. If $\Psi$ is $C^1$
we say that $\Psi$ is continuous if $\Psi_{M,N}f$ depends continuously on the metrics
of $M$ and $N$, for every pair $M$, $N$
(here we consider $\Psi f$ varying in the $C^1$ topology and the metrics varying in the $C^2$ topology). \\

$\Psi$ is 
cover-invariant if $\Psi\, \bar{f}=\overline{\Psi f}$
for every finite cover $\bar{f}: \bar{M}\ra\bar{N}$ of any $f:M\ra N$. Then we have:
\begin{enumerate}
\item[{a.}]  If $\Psi$ is  $C^1$,  then versions of Theorems 3 - 6 hold for $\Psi$-special maps.

\item[{b.}]  If $\Psi$ is a $C^1$, cover-invariant homotopy special correspondence,  then versions of 
Theorems 3 - 7 hold for $\Psi$-special maps.

\item[{c.}]  If, in addition,  $\Psi$ is continuous,  then versions of Theorems 3 - 10 hold for $\Psi$-special maps.
\end{enumerate}
\vspace{.6in}

{\large {\bf 4. Ricci flow and pinched negative curvature.}}\\

Until now we have dealt with processes that produce some special type
of map, e.g harmonic maps or natural maps. Now we discuss some processes that
produce a special type of metrics: Einstein metrics, that is, metrics of constant
Ricci curvature. As argued in the introduction of Besse's book 
``Einstein Manifolds" \cite{B}, Einstein metrics are ideal in the sense that
they are not as general as metrics of constant scalar curvature, and they
are not as restrictive as metrics of constant sectional curvatures.
Note that every metric of constant sectional curvature is an Einstein metric.
In particular every hyperbolic manifold is an Einstein manifold
(i.e a Riemannian manifold with a complete Einstein metric).
In dimension three ``constant Ricci curvature" is equivalent to
``constant sectional curvature"; hence every 3-dimensional Einstein
manifold is a space-form.\\

The most well known method for obtaining Einstein metrics is the Ricci flow method
introduced by Hamilton in his seminal paper \cite{H}. 
Starting with an arbitrary smooth Riemannian metric $h$ on a closed smooth $n$-dimensional 
manifold $M^n$, he considered the evolution
equation $$ \frac{\partial}{\partial t}\, h=\frac{2}{n}\, r \, h \, -\, Ric$$

\noindent where $r =\int R\, d\mu /\,  \int d\mu$ is the average scalar curvature
($R$ is the scalar curvature) and $Ric$ is the Ricci curvature tensor of $h$.
Hamilton then spectacularly illustrated the success of this method by proving,
when $n=3$, that if the initial Riemannian metric has strictly positive
Ricci curvature it evolves through time to a positively curved Einstein
metric $h_\infty$ on $M^3$. And, because $n=3$, $(M^3,h_\infty )$ is a spherical
space-form; i.e. its universal cover is the round sphere. Following Hamilton's approach
G. Huisken \cite{Hu}, C. Margerin \cite{Ma} and S. Nishikawa \cite{N} proved that, for every $n$,
Riemannian $n$-manifolds whose sectional curvatures are pinched close to +1
(the pinching constant depending only on the dimension)  can be deformed, through the Ricci 
flow, to a spherical-space form.\\

Ten years after Hamilton's results appeared, R. Ye \cite{Ye} studied the Ricci flow when the initial
Riemannian metric $h$ is negatively curved and proved that a negatively curved
Einstein metric is strongly stable; that is, the Ricci flow starting near such a
Riemannian metric $h$ converges (in the $C^\infty$ topology)
to a Riemannian metric isometric to $h$, up to scaling.
(We introduce the notation $h\equiv h'$ for two Riemannian metrics that are isometric
up to scaling.) 
In \cite{Ye} R. Ye also proved that sufficiently pinched to -1 manifolds can be deformed,
through the Ricci flow, to hyperbolic manifolds, but the pinching constant in his Theorem
depends on other quantities (e.g the diameter or the volume).
Ye's paper was motivated by the problem
on whether the Ricci flow can be used to deform every sufficiently pinched
to -1 Riemannian metric to an Einstein metric
(the pinching constant depending only on the dimension). 
His paper partially implements a scheme
proposed by Min-Oo \cite{M}.\\

We say the the Ricci flow for a negatively
curved Riemannian metric $h$ {\it converges smoothly}  if the Ricci flow, starting at $h$, 
is defined for all $t$ and converges (in the $C^\infty$ topology) to a well defined 
negatively curved (Einstein) metric.
The next Theorem shows the existence of pinched negatively curved metrics for which the
Ricci flow does not converge smoothly. \\

\noindent {\bf Theorem 14. \cite{FO2}} {\it Given $n>10$ and $\epsilon >0$ there is a closed smooth
$n$-dimensional manifold $N$ such that

(i) $N$ admits a hyperbolic metric

(ii) $N$ admits a Riemannian metric $h$ with sectional curvatures in
$[-1-\epsilon ,-1+\epsilon]$ for which the Ricci flow does not converge smoothly.}\\

The key ingredients in the proof of Theorem 14 is Theorem 11 and the fact that
the Ricci flow satisfies the following properties:
\begin{enumerate}
\item[{1.}] Hyperbolic metrics are fixed by the Ricci flow.

\item[{2.}] The Ricci flow preserves isometries.

\item[{3.}] The limit of the Ricci flow (in case it exists) is cover invariant.

\item[{4.}] The Ricci flow depends continuously on initial data.
\end{enumerate}

\noindent {\bf Remarks.}

\noindent {\bf 1.}  By ``cover invariant" we mean the following:
let $g$ be a metric on $M$ and $p:\bar{M}\ra M$ a cover. If $g_t$ is
the Ricci flow starting at $g$ and converging to $g_{\infty}$, then
the Ricci flow starting at $p^*g$ converges to $p^*g_{\infty}$.

\noindent {\bf 2.} Property 4 does not state that the {\it limit} of the 
Ricci flow is continuous on initial data. \\

As we did in section 3, this can also be generalized. Before we give a definition
trying to formalize a general process for obtaining Einstein metrics 
on $\epsilon$-pinched to -1 Riemannian manifolds, we establish
some notation. ${\cal M}_P$ will denote the space of all Riemannian metrics on a smooth manifold $P$.
For $\epsilon >0$, let ${\cal M}^{\epsilon}_P$ denote the space of $\epsilon$-pinched to -1
Riemannian metrics on $P$. Also, ${\cal E}_P\subset {\cal M}_P$ will denote the space of negatively curved 
Einstein metrics on $P$. Recall that ${\cal E}_P/\equiv$ is discrete, see \cite{B}, p.357.\\

\noindent {\bf Definition.} Let $\epsilon > 0$ and $n$ be a positive integer. An {\it Einstein correspondence} 
$\Phi :{\cal M}^\epsilon\rightarrow {\cal E}$ {\it for $n$-dimensional manifolds}  
is a family of maps $\Phi_P :{\cal M}^\epsilon_P\rightarrow{\cal E}_P$, for every $n$-dimensional 
manifold $P$ for which  ${\cal M}^\epsilon_P$ is not empty.
We say that $\Phi$ is {\it cover-invariant} if $\Phi ( p^*g)\, =\, p^*(\Phi (g))$ for every finite cover
$p: P\rightarrow Q$ and $g\in{\cal M}^\epsilon_Q$, for which $\Phi_Q$ is defined.\\

We say that $\Phi$ is {\it continuous} if each $\Phi_P :{\cal M}^\epsilon_P\rightarrow{\cal E}_P$
is continuous. Here we consider ${\cal M}^\epsilon_P$ with the $C^\infty$ topology and ${\cal E}_P$
with the $C^2$ topology.\\

Let $h$, $h'$ $\in {\cal M}_P$. Write $h\equiv_0h'$ provided $(P,h)$ is isometric to $(P,h')$, up to scaling,
via an isometry homotopic to $id_P$. Notice that the fibers of ${\cal E}_P/\equiv_0\,\, \rightarrow \,\, {\cal E}_P/\equiv$
are discrete; and hence ${\cal E}_P/\equiv_0$ is also discrete.\\

\noindent {\bf Theorem 15.  \cite{FO2}} {\it 
Suppose that there are  $\epsilon >0$ and $n>10$ for which there exists a cover-invariant
Einstein correspondence $\Phi$. Then there is a closed
$n$-dimensional Riemannian manifold $N$, with metric $h\in{\cal M}^\epsilon_N$, for which 
the Einstein metric $\Phi ( h)$ is unreachable by the Ricci flow starting at $h$.}\\

The proof of Theorem 15 is similar to the proof of Theorem 14, see \cite{FO2}.\\

\noindent {\bf Theorem 16. \cite{FO2}} {\it 
Suppose that there are  $\epsilon >0$ and $n\geq 6$ for which there exists an
Einstein correspondence $\Phi$. Then there is a closed
$n$-dimensional manifold $N$ that admits, at least, two non-isometric (even after scaling)
negatively curved Einstein metrics. Moreover, one metric can be chosen to be 
hyperbolic.}\\

This Theorem is easily deduced from Theorem 8. We reproduce the proof:\\

\noindent {\bf Proof.} 
From Theorem 8 we have the following.\\

There are closed connected smooth manifolds $M_0$, $M_1$, $N$,
of dimension $n$, Riemannian metrics $g_0$, $g_1$ on $M_0$
and $M_1$, respectively, and smooth two-sheeted covers 
$p_0 : N\rightarrow M_0$,   $p_1 : N\rightarrow M_1$ such that: \\

(1) $M_0$ and $M_1$ are homeomorphic but not $PL$-homeomorphic.

(2) $g_0$ is hyperbolic

(3) $g_1$ has sectional curvatures in $[-1-\epsilon ,-1+\epsilon]$.\\

Now, note that the metrics $g_1$ and $\Phi(g_1)$ are not hyperbolic, otherwise, by Mostow's Rigidity
Theorem, $M_0$ would be diffeomorphic to $M_1$, which contradicts (1) above. Hence  $p_1^* (\Phi (g_{1}) )$
is not hyperbolic either, while $p_0^*(g_0)$ is hyperbolic.
Then the two non-isometric negatively curved Einstein metrics on $N$ are
$p_0^*(g_0)$ and $p_1^* (\Phi (g_{1}) )$. This proves Theorem 16.\\

The general form of the following Theorem was suggested to us by Rugang Ye.\\

\noindent {\bf Theorem 17. \cite{FO2}} {\it A cover-invariant Einstein correspondence cannot be
continuous.}\\

\noindent {\bf Remark.} Note that we are not assuming that $\Phi$ fixes hyperbolic metrics.
If we assumed that $\Phi (hyperbolic\,\,\, metric)=(hyperbolic\,\,\, metric),$ the Theorem 
then would be easily deduced as before.\\

Since Ricci flow and elliptic deformation are cover-invariant continuous (analytic) processes, it follows
from Theorem 17 that they cannot be used, at least directly, to find Einstein metrics on
$\epsilon$-pinched to -1 Riemannian manifolds.\\

{\bf Dedication.} This article is respectfully dedicated to the memory of Armand Borel
whose conjecture that a closed aspherical  manifold is determined (up to homeomorphism)
by its fundamental group was one motivation for the research surveyed here.

F.T. Farrell

SUNY, Binghamton, N.Y., 13902, U.S.A.\\

P. Ontaneda

UFPE, Recife, PE 50670-901, Brazil


\begin{thebibliography}{99}
\bibitem{A'} S. I. Al'ber, {\em Spaces of mappings into manifold of negative
curvature}, Dokl. Akad. Nauk USSR {\bf168}, (1968) 13-16.



\bibitem{B}  A. L. Besse, {\em Einstein Manifolds}, Ergebnisse Series vol. {\bf 10},
Springer-Verlag, Berlin, 1987.

\bibitem{BCG} G. Besson, G. Courtois and S. Gallot, {\em Minimal entropy and Mostow's rigidity
Theorems}, Ergodic Theory \& Dynam. Sys. {\bf 16} (1996), 623-649.


   
\bibitem{BL} D. Burghelea and R. Lashof, {\em Stability of 
concordances and suspension homeomorphism}, Ann. of Math. (2) {\bf 
105} (1977), 449-472.

\bibitem{BK} K. Burns and A. Katok, {\em Manifolds with non-positive curvature}, Ergodic
Theory \& Dynam. Sys. {\bf 5} (1985), 307-317.


%\bibitem{Bor63} A. Borel, {\em Compact Clifford-Klein forms of symmetric spaces}, Topology {\bf2}, (1963) 111-122.



%\bibitem{CR} P.E. Conner and F. Raymond, {\em Deforming homotopy equivalences to homeomorphisms in aspherical
%manifolds}, Bull. AMS, {\bf 83} (1977) 36-85. 


\bibitem{Cor} K. Corlette, {\em Archimedean superrigidity and hyperbolic geometry},
Ann. of Math. {\bf 135} (1992), 165-182.


\bibitem{E} R.D. Edwards, {\em The topology of manifolds and cell-like maps}, in Proc. of the
ICM (Helsinki, 1978), pp. 111-127, Acad. Sci. Fennica, Helsinki, 1980.



\bibitem{EL1} J. Eells and L. Lemaire, Selected topics in harmonic maps, CBMS
Regional Conf. Series 50, Amer. Math. Soc., Providence, R.I\., 1983.


\bibitem{EL2} J. Eells and L. Lemaire, {\em Another report on harmonic maps}, Bull. of LMS,
{\bf 20} (1988), 385-524.


\bibitem{ES} J. Eells and J. H. Sampson, {\em Harmonic mappings of Riemannian manifolds}, Amer. J. Math. {\bf86}, (1964) 109-160.


\bibitem{FJ} F.T. Farrell and L.E. Jones, {\em K-Theory and dynamics II}, Ann. of Math.
{\bf 126} (1987) 451-493.



\bibitem{FJ1} F.T. Farrell and L.E. Jones, {\em Negatively curved manifolds with exotic smooth structures}, J. Amer. Math. Soc. 
{\bf 2} (1989) 899-908.

\bibitem{FJ2} F.T. Farrell and L.E. Jones, {\em Rigidity in geometry and
topology}, Proc. of the International Congress of Mathematicians,
Vol. I, II (Kyoto, 1990), Math. Soc. Japan, Tokyo (1991)   653-663.


\bibitem{FJ3} F.T. Farrell and L.E. Jones, {\em Some non-homeomorphic harmonic homotopy equivalences}, Bull. London Math. Soc. 
{\bf 28} (1996), 177-180.




\bibitem{FJO1} F.T. Farrell, L.E. Jones and P. Ontaneda, {\em Hyperbolic manifolds with negatively curved exotic 
triangulations in dimension larger than five}. Jour. Diff. Geom. {\bf 48} (1998) 319-322.

\bibitem{FJO2} F.T. Farrell, L.E. Jones and P. Ontaneda, {\em Examples of non-homeomorphic
harmonic maps between negatively curved manifolds}, Bull. London Math. Soc. {\bf 30} 
(1998) 295-296.

\bibitem{FOR} F.T. Farrell, P. Ontaneda and M.S. Raghunathan, {\em Non-univalent harmonic maps homotopic
to diffeomorphisms}, Jour. Diff. Geom. {\bf 54} (2000) 227-253.

\bibitem{FO1} F.T. Farrell and P. Ontaneda, {\em Cellular harmonic maps which are not diffeomorphisms}.
Submitted for publication. arXiv:org/mathDG/0311175.

\bibitem{FO2} F.T. Farrell and P. Ontaneda, {\em A caveat on the convergence of the Ricci flow
for negatively curved manifolds},
Submitted for publication. arXiv:org/mathDG/0311176.

\bibitem{GS} M. Gromov and R. Schoen, {\em Harmonic maps into singular spaces and $p$-adic
superrigidity of lattices in groups of rank one}, Inst. Hautes {\'E}tudes Sci. Publ. Math. 
{\bf 76} (1992) 165-246.


\bibitem{H} R. Hamilton, {\em Three-manifolds with positive Ricci curvature}, 
Jour. Diff. Geom. {\bf  17} (1982) 255-306.



\bibitem{Har} P. Hartman, {\em On homotopic harmonic maps}, Canad. J. Math. {\bf19}, (1967) 673-687.


\bibitem{Hat} A.E. Hatcher, {\em Concordance spaces, higher simple homotopy theory, and
applications}, Proc. Symp. Pure Math. {\bf 32} (1978), 3-21.



\bibitem{Her} L. Hern{\'a}ndez, {\em K{\"a}hler manifolds and 1/4-pinching}, Duke
Math. J. {\bf 62} (1991) 601-611.



\bibitem{Hu} G. Huisken, {\em Ricci deformation of a metric on a Riemannian manifold}, Jour. Diff.Geom.
{\bf 21} (1985) 47-62.


\bibitem{I} K. Igusa, {\em Stability Theorems for pseudoisotopies}, 
K-theory {\bf 2} (1988) 1-355.


\bibitem{JY} J. Jost and S.-T. Yau, {\em Harmonic maps and superrigidity}, Proc.
Sympos. Pure Math. {\bf 54} Amer. Math. Soc., Providence, R.I., (1993)
245-280.


\bibitem{KPW} A. Katok, M. Pollicott and H. Weiss, {\em Differentiability of entropy for Anosov geodesic flows},
Bull. Amer. Math. Soc. {\bf 22} (1990) 285-293.

\bibitem{KiSi} R.C. Kirby and L.C. Siebenmann, {\em Foundational Essays on 
Topological Manifolds, Smoothings, and Triangulations}, Annals of Math. Studies
no.{\bf88}, Princeton University Press, Princeton (1977).

%\bibitem{LawYa} H.B. Lawson and S.T. Yau, {\em Compact manifolds of nonpositive
%curvature}, J. Diff. Geom. {\bf7}, (1972) 211-228.


\bibitem{Ma} C. Margerin, {\em  Pointwise pinched manifolds are space forms}, Proc. Sympos. Pure Math.
vol. {\bf 44}, Amer. Math. Soc., Providence, R.I., 1986, 307-328.

\bibitem{M} M. Min-Oo, {\em Almost Einstein manifolds of negative Ricci curvature}, Jour. Diff. Geom.
{\bf 32} (1990) 457-472.

\bibitem{Mill} J.J. Millson, {\em On the first Betti number of a constant 
negatively curved manifold}, Ann. of Math. {\bf 104} (1976) 235-247.
\bibitem{MR} J.J. Millson and M.S. Raghunathan, {\em Geometric construction
of cohomology for arithmetic groups I}, Proc. Indian Acad. Sci., Vol.{\bf90}, no.{\bf2} (1981) 103-123.


\bibitem{MSY} N. Mok, Y.-T. Siu and S.-K. Yeung, {\em Geometric superrigidity},
Invent. Math. {\em 113} (1993), 57-83.


\bibitem{Mos} G.D. Mostow, {\em Quasi-conformal mappings in $n$-space    
and the rigidity of hyperbolic space forms}, Inst. Hautes \'Etudes
Sci. Publ. Math. {\bf 34} (1967), 53-104.



%\bibitem{MoTa} G.D. Mostow and T. Tamagawa, {\em On the compactness of arithmetically defined homogeneous spaces}, Ann. of Math. {\bf76} (1962) 463-466.\bibitem{Yau} S. T. Yau, {\em Seminar on differential geometry}, Ann. of Math. Stud., {\bf102}, Princeton Univ. Press, Princeton, NJ, 1982.

\bibitem{N} S. Nishikawa, {\em Deformation of Riemannian metrics a manifolds with bounded curvature
ratios}, Proc. Sympos. Pure Math. vol. {\bf 44}, Amer. Math. Soc., Providence, R.I., 1986, 343-352.



\bibitem{O} P. Ontaneda, {\em Hyperbolic manifolds with negatively curved exotic triangulations in dimension six}, 
J. Diff. Geom. {\bf 40} (1994), 7-22. 


\bibitem{Sa} J. Sampson, {\em Some properties and applications of harmonic
mappings}, Ann. Scient. Ec. Norm. Sup. {\bf 11} (1978), 211-228.


\bibitem{Sa2} J. Sampson, {\em Applications of harmonic
maps to K{\"a}hler geometry}, Cont. Math. {\bf 49} (1986), 125-133.




\bibitem{Siu} Y.-T. Siu, {\em The complex-analyticity of harmonic maps and the
strong rigidity of compact K{\"a}hler manifolds}, Ann. of Math. {\bf 112} (1980),
73-111.


\bibitem{SC1} M. Scharlemann and L. Siebenmann, {\em The Hauptvermutung for smooth singular homeomorphisms}, 
in Manifolds Tokyo 1973, Akio Hattori ed. Univ. of Tokyo Press, 85-91.



\bibitem{SC2} M. Scharlemann, {Smooth CE maps and smooth homeomorphisms}, in Algebraic and Geometric Topology, 
A. Dold and B. Eckmann eds. Lecture Notes in Mathematics n.664,
234-240.


\bibitem{SY} R. Schoen and S.-T. Yau, {\em On univalent harmonic maps between
surfaces}, Inv. Math. {\bf 44} (1978), 265-278.


\bibitem{Sh} S.A. Shcherbakov, {\em On the degree of smoothness of horospheres, radial fields and horospherical
coordinates on a Hadamard manifold}, Soviet. Math. Dokl. {\bf 28} (1983) 233-237.

\bibitem{Si} L.C. Siebenmann, {\em Approximating cellular maps by homeomorphisms}, Topology 
{\bf 11} (1973), 271-294.



\bibitem{Sta} C.W. Stark, {\em Surgery theory and infinite fundamental groups}, 
Ann. of Math. Studies {\bf 145}, Vol.1, 239-252.


\bibitem{V} M. Varisco, personal communication.

\bibitem{Yau} S.-T. Yau, {\em Seminar on differential geometry}, Ann.of Math. Stud., {\bf 102},
Princeton Univ. Press, Princeton, NJ, 1982. 




\bibitem{YZ} S.-T. Yau and F. Zheng, {\em  Negatively 1/4-pinched Riemannian
metric on a compact K{\"a}hler manifold}, Invent. Math. {\bf 103} (1991),
527-535.

\bibitem{Ye} R. Ye, {\em Ricci flow, Einstein metrics and space forms}, 
Transactions of the AMS {\bf 338} (1993) 871-896. \\





\end{thebibliography}
\end{document}